\documentclass[twocolumn,showpacs,%
  nofootinbib,aps,superscriptaddress,%
  eqsecnum,prd,notitlepage,showkeys,10pt]{revtex4-1}

\usepackage{amssymb}
\usepackage{amsmath}
\usepackage{graphicx}
\usepackage{dcolumn}
\usepackage{hyperref}
\usepackage{bbold}

\newtheorem{theorem}{Theorem}[section]

\begin{document}

\title{On the Fermat's Last Theorem and the Dirac equation}
\author{Mart\'in D. Arteaga Tupia}
\email{martin77@if.usp.br}
\affiliation{Instituto de F\'isica, Universidade de S\~ao Paulo, 05314-970 S\~ao Paulo, Brazil}

\affiliation{Callao Theoretical Physics group, Callao, Per\'u}

\begin{abstract}
In the present paper we study, in a mathematically non-formal way, the validity of the Fermat's Last Theorem (FLT) by generalizing the usual procedure of extracting the square root of non convenient objects initially introduced by P. A. M. Dirac in the study of the linear relativistic wave equation.   
\end{abstract}

\maketitle
%
%
%
\section{Introduction}
Most people know that around 1637 Fermat wrote the following statement \cite{Panchi}
\begin{quote}
\it``It is impossible to separate a cube into two cubes, or a fourth power into two fourth powers, or in general, any power higher than the second, into two like powers. I have discovered a truly marvellous proof of this, which this margin is too narrow to contain.''
\end{quote}
 in the margin of the Arithmetica of Diophantus. Such a statement in modern language is translated to

\begin{theorem}
With $n, x, y, z \in \mathbf{N}$ and $n > 2$ the equation $x^n + y^n = z^n$ has no solutions.
\end{theorem}

In a monumental (and very large) work, during 1994, A. Wiles found a proof for the above theorem using many modern techniques that are not easy to digest for the non-mathematician community. 

On the other hand, in 1928, P. A. M. Dirac by looking for a linear relativistic wave equation version of the Klein-Gordon one found the currently best known Dirac equation,
\begin{eqnarray}\label{1}
(i\partial\!\!\!/ - m ) \psi = 0.
\end{eqnarray}
In the road to find the above linear version of $(\square - m^2)\phi = 0$, one usually introduces the {\it Gamma matrices} and their algebra
\begin{eqnarray}\label{2}
\{\gamma^\mu, \gamma^\nu\}=2\eta^{\mu\nu}.
\end{eqnarray}
However, Dirac used another notation, based on the matrix coefficients $\hat{\alpha}^i$ and $\hat{\beta}$ \cite{Dirac}. Such objects was introduced in order to write a linear version of the square root of the hamiltonian operator
\begin{eqnarray}\label{3}
\hat{\mathcal{H}} = \hat{\alpha}^i\, p^i + \hat{\beta}\, m
\end{eqnarray}
and by imposing that the eingenvalues of $\hat{\mathcal{H}}$ and $p^i $ satisfy $E^2=p^2+m^2$, it is found 
\begin{eqnarray}
(\hat{\alpha}^i)^2 &=& \mathbb{1},\label{4}\\ 
(\hat{\beta})^2 &=& \mathbb{1},\label{5} \\
\hat{\alpha}^i \hat{\beta} + \hat{\beta} \hat{\alpha}^i &=&0.\label{6}
\end{eqnarray}

In the following sections we will use an analogous procedure to show that the equation $x^n + y^n = z^n$, with $n, x, y, z \in \mathbf{N}$, can not be solved for $n>2$, because we can not have a consistent and solvable set of equations as eqs.(\ref{4}), (\ref{5}) and (\ref{6}).

\section{Linearizing $z^n = x^n + y^n$}

Let us start by writing the following linearized {\it ansatz}  
\begin{eqnarray}\label{2.1}
z\,\Gamma^z_{n} = x \, \Gamma^x_{n} + y \, \Gamma^y_{n}
\end{eqnarray}
where $\Gamma^i_{n}$ , $i= x, y, z$, are matrices. Therefore, if $z^n = x^n + y^n$ is transformed into eq.(\ref{2.1}) with well defined $\Gamma^i_{n}$, it is easy to verify that the above equation admit three integers $x, y$ and $z$ simultaneously as a solution.

Therefore, the problem is now transformed into:
\begin{quote}
Is it possible to find well defined $\Gamma^i_{n}$ for arbitrary $n$?
\end{quote}
\subsection*{Case n=1}
This is the trivial case, $\Gamma^i_{1} = \mathbb{1}$ for $i= x, y, z$. 

\subsection*{Case n=2}
This Pythagorean case is given by
\begin{eqnarray}\label{2.2}
z^2 = x^2 + y^2,
\end{eqnarray}
we then write 
\begin{eqnarray}\label{2.3}
z\,\Gamma^z_{2} = x \, \Gamma^x_{2} + y \, \Gamma^y_{2}
\end{eqnarray}
and in order it will satisfy $z^2 = x^2 + y^2$ we found
\begin{eqnarray}
(\Gamma^z_{2})^2=\mathbb{1}, (\Gamma^x_{2})^2=\mathbb{1}, (\Gamma^y_{2})^2 &=&\mathbb{1},\label{2.4}\\
\Gamma^x_{2}\Gamma^y_{2} + \Gamma^y_{2} \Gamma^x_{2} = 0\label{2.5}.
\end{eqnarray}
We note from eq.(\ref{2.4}) that we have three unknown objects $\Gamma^z_{2}, \Gamma^x_{2}$ and $\Gamma^y_{2}$, and also three equations eqs.(\ref{2.2}), (\ref{2.3}) and (\ref{2.5}), which relate them. Therefore, this is a compatible set of equations.

To find the matrices explicitly we fix $\Gamma^x_{2}$ to be
\begin{eqnarray}\label{2.6}
\Gamma^x_{2} =
\begin{pmatrix}
 0 & 1 &\\
 1 & 0 &
\end{pmatrix}
\end{eqnarray}
which is a Pauli matrix, we also choose that matrix because $(\Gamma^x_{2})^2=\mathbb{1}$ and has zero trace. By using eq.(\ref{2.5}) we can get another Pauli matrix, say
\begin{eqnarray}\label{2.7}
\Gamma^y_{2} =
\begin{pmatrix}
 1 & 0 &\\
 0 & -1 &
\end{pmatrix}.
\end{eqnarray}
Using eq.(\ref{2.3}), we are able to obtain 
\begin{eqnarray}\label{2.8}
\Gamma^z_{2} =
\begin{pmatrix}
 y/z & x/z &\\
 x/z & -y/z &
\end{pmatrix}
\end{eqnarray}
and it is easy to check that $(\Gamma^z_{2})^2=\mathbb{1}$, as expected. For instance, in the case of the well known $3^2 + 4^2 = 5^2$, such a matrix is given by 
\begin{eqnarray}\label{2.9}
\Gamma^z_{2} =
\begin{pmatrix}
 4/5 & 3/5 &\\
 3/5 & -4/5 &
\end{pmatrix}.
\end{eqnarray}
\subsection*{Case n=3}
After linearizing
\begin{eqnarray}\label{2.10}
z^3 = x^3 + y^3,
\end{eqnarray}
we have to obtain
\begin{eqnarray}\label{2.11}
z\,\Gamma^z_{3} = x\,\Gamma^x_{3} + y\,\Gamma^y_{3}
\end{eqnarray}
and after some algebra we arrive to the following relations
\begin{eqnarray}
(\Gamma^z_{3})^2=\mathbb{1}, (\Gamma^x_{3})^2=\mathbb{1}, (\Gamma^y_{3})^2 &=&\mathbb{1},\label{2.12}\\
(\Gamma^x_{3})^{2}\Gamma^y_{3} +  \Gamma^x_{3} \Gamma^y_{3} \Gamma^x_{3} +\Gamma^y_{3} (\Gamma^x_{3})^{2} &=& 0,\label{2.13}\\
(\Gamma^y_{3})^{2}\Gamma^x_{3} +  \Gamma^y_{3} \Gamma^x_{3} \Gamma^y_{3} +\Gamma^x_{3} (\Gamma^y_{3})^{2} &=& 0.\label{2.14}
\end{eqnarray}
As before, we have three unknown variables, the matrices $\Gamma^x_{3}, \Gamma^y_{3}$ and $\Gamma^z_{3}$. Nevertheless, now the number of equations are four: eqs.(\ref{2.10}), (\ref{2.11}), (\ref{2.13}) and (\ref{2.14}). Consequently, such a system is incompatible.

In this way, we conclude that this case fail to satisfy eq.(\ref{2.10}) for $x, y, z \in \mathbf{N}$.

For $n> 3$ the situation is worst, the number of equations grows while the number of unknown variables remains fixed to be three: $\Gamma^z_{n}, \Gamma^z_{n}$ and $\Gamma^z_{n}$.

\subsection*{General case}

In general, we want to have the certainty on what is the maximum value of $n$, which allows to have three integer values of $x, y$ and $z$ at the same time, such that $x^n +y^n = z^n$.

For $n=1$, this is trivial because it is just a sum. Therefore, if we linearize $x^n +y^n = z^n$ in some way, to find $x, y$ and $z$ integers is easy. In this way, we have to figure out under what conditions, for $n$, we can linearize the above equation.

Our main equation to be solved is thus
\begin{eqnarray}\label{2.15}
z^n = x^n +y^n ,
\end{eqnarray}
the second one is 
\begin{eqnarray}\label{2.16}
\Gamma^z_n\, z = \Gamma^x_n\, x + \Gamma^y_n\, y.
\end{eqnarray}
After powering to $n$ eq.(\ref{2.16}), we arrive to 
\begin{eqnarray}
(\Gamma^z_n\, z)^n &=& (\Gamma^x_n\, x + \Gamma^y_n\, y)^n\nonumber\\
&=& (\Gamma^x_n)^n \, x^n +\sum_{k=1}^{n-1}Per[(\Gamma^x_n)^{n-k}\,(\Gamma^y_n)^k]x^{n-k}\,y^k\nonumber\\ 
&&+ (\Gamma^y_n)^n \, y^n,\label{2.17}
\end{eqnarray}
where we have defined $\sum_{p,q}\,Per(\hat{\alpha}^p\, \hat{\beta}^{q})$ as the sum of the all possible $p+q$ permutations of the non commuting objects $\hat{\alpha}$ and $\hat{\beta}$.

As such, we have to have
\begin{eqnarray}\label{2.18}
\sum_{\text{fixed}\,\,k}Per[(\Gamma^x_n)^{n-k}\,(\Gamma^y_n)^k] = 0,
\end{eqnarray}
this is a set of $n-1$ equations, and together with eqs.(\ref{2.15}), (\ref{2.16}) form a total of $n+1$ conditions. Hence, recalling that the number of the unknown variables are the three matrices $\Gamma_{n}^{i}$, we have that the system will be compatible if
\begin{eqnarray}
n+1 \leq 3 \rightarrow n \leq 2 .
\end{eqnarray}

Therefore, we conclude that the equation $z^n = x^n + y^n$ can admit $n,x,y, z \in \mathbf{N}$, as solutions, only for $n<3$.

\section{Final Comments}

Because of the algebra, this paper does not fit in the margin of a book.
  


\end{document}